\documentclass[aos,preprint]{imsart}

\usepackage{amsmath, amsthm, amssymb}
\usepackage{graphicx}
\usepackage{verbatim}
\usepackage{caption}
\usepackage{subcaption}
\usepackage{fancyvrb}
\usepackage{enumerate}

\RequirePackage[numbers]{natbib}
\RequirePackage[colorlinks,citecolor=blue,urlcolor=blue]{hyperref}

\usepackage[margin=1.5in]{geometry}

\theoremstyle{plain}
\newtheorem{prop}{Proposition}

\newtheorem{lemm}[prop]{Lemma}
\newtheorem{theo}[prop]{Theorem}

\theoremstyle{definition}
\newtheorem{exam}{Example}
\newtheorem{defi}[exam]{Definition}

\theoremstyle{remark}

\usepackage{StefanTexCommandsV2}
\graphicspath{{./figures/}}

\newcommand{\rf}{\operatorname{RF}}
\newcommand{\wrf}{\proj{\rf}}

\newcommand{\hVIJ}{\widehat{V}_{IJ}}

\newcommand{\proj}[1]{\mathring{#1}}

\newcommand{\figw}{0.8\textwidth}

\newcommand{\SPACEBIG}{1}
\newcommand{\SPACESMALL}{1}

\def\spacingset#1{\renewcommand{\baselinestretch}%
{#1}\small\normalsize} \spacingset{1}

\begin{document}

\begin{frontmatter}

\title{Asymptotic Theory for Random Forests}
\runtitle{Asymptotic Theory for Random Forests}

\begin{aug}
\author{\fnms{Stefan} \snm{Wager}\corref{}\ead[label=e1]{swager@stanford.edu}\thanksref{t1}}
\runauthor{S. Wager}
\thankstext{t1}{I am grateful for many helpful conversations with Brad Efron, Trevor Hastie, and Guenther Walther. This work was supported by a B.C. and E.J. Eaves Stanford Graduate Fellowship.} 
\address{Department of Statistics \\ Stanford University \\ Stanford, CA-94305 \\ \printead{e1}}
\affiliation{Stanford University}
\end{aug}


\begin{abstract}
Random forests have proven to be reliable predictive algorithms in many application areas. Not much is known, however, about the statistical properties of random forests. Several authors have established conditions under which their predictions are consistent, but these results do not provide practical estimates of random forest errors. In this paper, we analyze a random forest model based on subsampling, and show that random forest predictions are asymptotically normal provided that the subsample size $s$ scales as $s(n)/n = o(\log(n)^{-d})$, where $n$ is the number of training examples and $d$ is the number of features. Moreover, we show that the asymptotic variance can consistently be estimated using an infinitesimal jackknife for bagged ensembles recently proposed by Efron (2014). In other words, our results let us both characterize and estimate the error-distribution of random forest predictions, thus taking a step towards making random forests tools for statistical inference instead of just black-box predictive algorithms.
\end{abstract}

\end{frontmatter}

\spacingset{\SPACEBIG}

{\bf Remark.} This manuscript is superseded by the paper `` Estimation and Inference of Heterogeneous Treatment Effects using Random Forests'' by Wager and Athey, available at \url{http://arxiv.org/pdf/1510.04342.pdf}. The new paper by Wager and Athey extends the asymptotic theory developed here, and applies it to causal inference in the potential outcomes framework with unconfoundedness. The present manuscript will remain online for archival purposes; however, for all proofs, please see the paper by Wager and Athey.

\section{Introduction}

Random forests, introduced by \citet{breiman2001random}, have become one of the most popular out-of-the-box prediction tools for machine learning. But despite the ubiquity of applications, the statistical properties of random forests are not yet fully understood. Simple questions such as ``What is the limiting distribution of random forest predictions as the number of training examples goes to infinity?'' and ``How can we get consistent estimates of the sampling variance of random forest predictions?'' still remain largely open. This paper provides answers to these questions for a large class of random forest models.

We study random forests with subsampling, defined as follows. Suppose that we have training examples $Z_i = \p{X_i, \, Y_i}$ for $i = 1, \, ..., \, n$, a test point $x$, and a potentially randomized regression tree predictor $T$ which makes predictions $\hy = T\p{x; \, Z_1, \, ..., \, Z_n}$. We can then turn this tree $T$ into a random forest by averaging it over $B$ random samples:
\begin{equation}
\label{eq:rfm}
\rf_s^B\p{x; \, Z_1, \, ..., \, Z_n} = \frac{1}{B} \sum_{b = 1}^B  T\p{x; \, Z^*_{b1}, \, ..., \, Z^*_{bs}} \text{ for some } s \leq n,
\end{equation}
where $\{Z^*_{b1}, \, ..., \, Z^*_{bs}\}$ form a uniformly drawn random subset of $\{Z_1, \, ..., \, Z_n\}$. As explained heuristically by \citet{breiman2001random} and more formally by, e.g., \citet{buhlmann2002analyzing}, the resampling performed by random forests improves single trees by smoothing their decision thresholds.

Our main result is that, for a large class of base learners $T$, the predictions made by a random forest of the form \eqref{eq:rfm} are asymptotically normal provided that $s\p{n}/n = o\p{\log(n)^{-d}}$, where $d$ is the dimension of the feature space. Moreover, under these conditions, we show that the asymptotic variance of the ensemble can be consistently estimated using a simple formula proposed by \citet{efron2013estimation}, the infinitesimal jackknife for bagged ensembles. Thus, our results let us take a step towards bringing random forest predictions into the realm of classical statistical inference.

\subsection{Related Work}

Breiman originally described random forests in terms of bootstrap sampling:
\begin{equation}
\label{eq:rf_plain}
\rf^B\p{x; \, Z_1, \, ..., \, Z_n} = \frac{1}{B} \sum_{b = 1}^B  T\p{x; \, Z^*_{b1}, \, ..., \, Z^*_{bn}},
\end{equation}
where the $Z^*_{bi} \simiid \{Z_1, \, ..., \, Z_n\}$ form a bootstrap sample of the training data. Random forests of the form \eqref{eq:rf_plain}, however, have proven to be remarkably resistant to classical statistical analysis. As observed by \citet{buja2006observations}, \citet{chen2003effects}, \citet{friedman2007bagging} and others, estimators of this form can exhibit surprising properties even in simple situations. In this paper, we avoid these pitfalls by basing our analysis on subsampling rather than bootstrap sampling.

Subsampling has often been found to be more easily amenable to theoretical analysis than bootstrapping \citep[e.g.,][]{politis1999subsampling}. Although Breiman originally built random forests using bootstrap aggregation, or bagging \citep{breiman1996bagging}, subsampled random forests have also been studied and found to work well \citep[e.g.,][]{buhlmann2002analyzing,strobl2007bias}.  Our present analysis is related to work by \citet{hall2005properties}, who study the properties of subsampled $k$-nearest neighbors predictors and show that they are consistent provided the subsample size $s\p{n}$ grows to infinity slower than $n$.

Most existing theoretical results about random forests aim to establish the consistency of random forest predictions \citep{biau2012analysis,biau2008consistency,breiman2004consistency,meinshausen2006quantile}. There has been much less work, however, on understanding the sampling variance of random forests.
One notable exception is a paper by \citet{lin2006random} that uses an adaptive nearest neighbors analysis of random forests to derive lower bounds for their variance. They then use this result to show that predictors of the form \eqref{eq:rf_plain} converge slowly in $n$ if the leaf-size of the trees $T$ is small.

Our work builds on recent contributions by \citet{duan2011bootstrap}, \citet{sexton2009standard} and \citet{wager2014confidence}, who seek to move beyond a black-box analysis of random forests by providing variance estimates for their predictions based on ideas like the bootstrap and the jackknife. These papers, however, only gave heuristic and experimental support for their proposed approaches. Here, we provide rigorous results about the asymptotic distribution of random forest predictions.

Finally, we note recent independent work by \citet{mentch2014ensemble}, who establish asymptotic normality of generic subsampled ensembles in a small-subsample-size regime with $s(n)/\sqrt{n} \rightarrow 0$, assuming a non-degeneracy condition on the first-order effects of the base learner. In comparison with this work, our result holds under the substantially weaker condition $s(n) / \log(n)^d \rightarrow 0$, and we show that the variance of random forest predictions can be estimated using the infinitesimal jackknife \citep{efron2013estimation,wager2014confidence}.

\section{Statistics of Random Forest Predictions}

The goal of this paper is to show how random forests can be made amenable to classical statistical analysis. Let
$\hy = \rf_{s}\p{x; \, Z_1, \, ..., \, Z_n}$ be the prediction made at $x$ by a random forest, with $B$ large enough that Monte Carlo effects have disappeared
\begin{equation}
\label{eq:infty_B}
\rf_s\p{x; \, Z_1, \, ..., \, Z_n} = \lim_{B \rightarrow \infty} \rf_s^B\p{x; \, Z_1, \, ..., \, Z_n}.
\end{equation}
We show that, for a large class of random forests, $\hy$ is consistent as $n$ grows to infinity and that
\begin{equation}
\label{eq:heuristic}
\hy \approxdot \nn\p{\EE{\hy}, \, \sigma^2\p{\hy}},
\end{equation}
provided that the random forest is trained using a subsample size $s(n)$ that is smaller than $n/\log(n)^d$; here $d$ is the number of features, $\sigma(\hy)$ is the noise scale of the predictions, and $\nn$ is the standard normal density.  Moreover, we show that $\sigma^2(\hy)$ can consistently be estimated from data.

We begin by introducing the infinitesimal jackknife estimate for $\sigma^2(\hy)$ below, and show it in action on a real dataset. In Section \ref{sec:main}, we then describe rigorous conditions under which \eqref{eq:heuristic} holds. We give an overview of our proof in Section \ref{sec:theor}, and provide a simulation study for the accuracy of the infinitesimal jackknife variance estimate in Section \ref{sec:sim}.

\subsection{The Infinitesimal Jackknife for Random Forests}
\label{sec:hvij}

In order to estimate $\sigma^2(\hy)$, we use the infinitesimal jackknife (or non-parametric delta-method) estimator $\hVIJ$ for bagging introduced by \citet{efron2013estimation}. This estimator can be computed using a particularly simple formula:
\begin{equation}
\label{eq:hvij}
\hVIJ\p{x; \, Z_1, \, ..., \, Z_n}  = \sum_{i = 1}^n \Cov[*]{T\p{x; \, Z^*_{1}, \, ..., \, Z^*_{n}}, \, N^*_i},
\end{equation}
where $N^*_i$ is the number of times $Z_i$ appears in the subsample used by $T$ and the covariance is taken with respect to the resampling measure. This formula arises by applying the original infinitesimal jackknife idea of \citet{jaeckel1972infinitesimal} to the resampling distribution.

The estimator $\hVIJ$ was studied in the context of random forests by \citet{wager2014confidence}, who showed empirically that the method worked well for many problems of interest. \citet{wager2014confidence} also emphasized that, when using $\hVIJ$ in practice, it is important to account for Monte Carlo bias; in our case, the appropriate correction is given in Section \ref{sec:sim}.

Our analysis provides theoretical backing to these results, by showing that $\hVIJ$ is in fact a consistent estimate for the variance $\sigma^2(\hy)$ of random forest predictions. The earlier work on this topic \citep{efron2013estimation,wager2014confidence} had only motivated the estimator $\hVIJ$ by highlighting connections to classical statistical ideas, but did not establish any formal justification for it.

\subsubsection{A First Example}

\spacingset{\SPACESMALL}

\begin{figure}
\centering
\includegraphics[width = \figw]{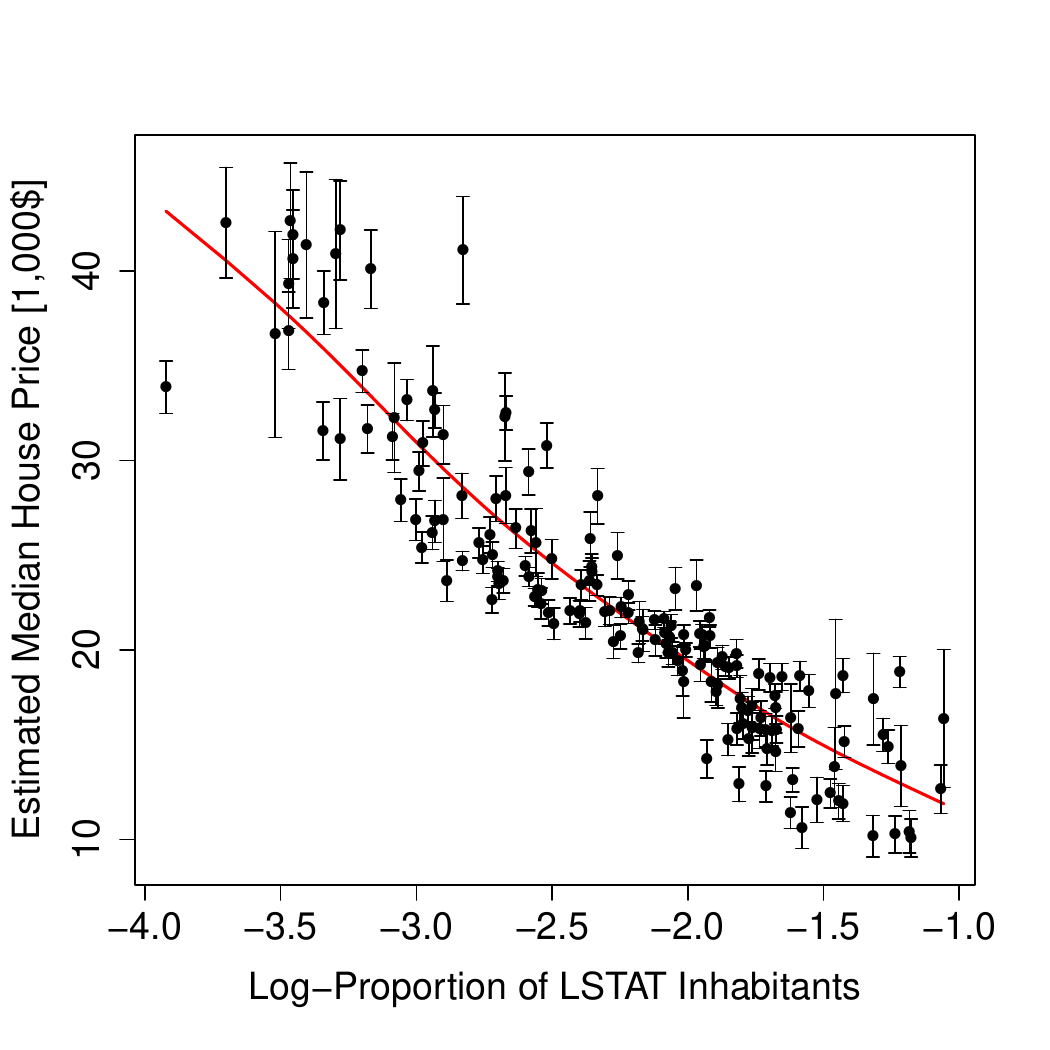}
\caption{Random forest predictions on the Boston housing data set \citep{harrison1978hedonic}, along with 95\% confidence intervals provided by the $\hVIJ$ estimator. The $x$-axis shows the log-proportion of lower status residents (LSTAT). The LSTAT proportion is the mean of the proportion of residents without some high school education and the proportion of male workers classified as laborers. We selected 337 out of 506 examples for training; the plot shows predictions made on the 169 remaining test examples. The random forest had a sub-sample size $s = 100$ and $B = 10,000$ replicates; otherwise, we used default settings for the \texttt{randomForest} package in \texttt{R} \citep{liaw2002classification}. The solid line is a smoothing spline with $df = 4$ degrees of freedom.}
\label{fig:lstat}
\end{figure}

\spacingset{\SPACEBIG}

To illustrate our inferential framework, we re-visit a classic regression data set: the Boston housing data set \citep{harrison1978hedonic}, with $d = 13$ features and $n = 337$ training examples. In Figure \ref{fig:lstat}, we plot random forest predictions for median house price against a measure of the proportion of lower status residents (LSTAT), along with 95\% confidence intervals provided by our theory.

We see that predictions with middle-range LSTAT values are all huddled near a smoothing spline drawn through the data, while the suburbs with extreme values of LSTAT have more scattered predictions. The error bars provided by $\hVIJ$ corroborate this observation: the size of the error bars roughly scales with the distance of the predictions from the smoothing spline trend line. Note that Figure \ref{fig:lstat} only shows 1 out of 13 predictors; this is why nearby points on the plot can have very different error bars.

\spacingset{\SPACESMALL}

\begin{figure}
\centering
\includegraphics[width = \figw]{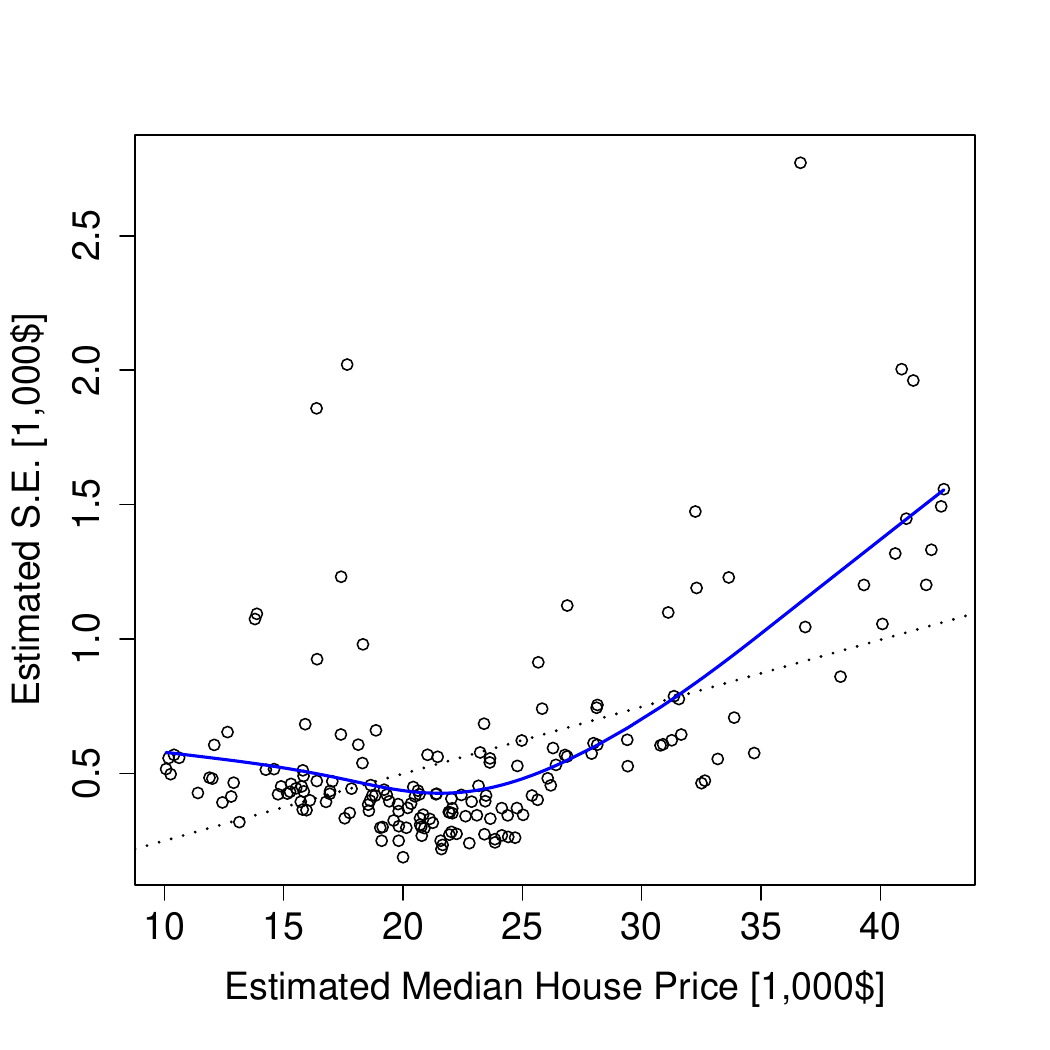}
\caption{Estimated standard errors for predictions on the Boston housing data set \citep{harrison1978hedonic}, as a function of the prediction itself. The solid line is a smoothing spline with $df = 4$ degrees of freedom, while the dotted line is the line connecting zero with the mean of the data. The experimental setup is described in the caption of Figure \ref{fig:lstat}.}
\label{fig:pred}
\end{figure}

\spacingset{\SPACEBIG}

Figure \ref{fig:pred} shows the relationship between random forest predictions $\hy$ and the estimated standard error $\hsigma(\hy)$. We might have expected for $\hsigma(\hy)$ to scale roughly proportionally to $\hy$. This behavior does in fact hold for large valuations $\hy$. However, at the low end of the range of $\hy$, we see that $\hsigma(\hy)$ starts climbing back up; in relative terms,  the coefficient of variation grows from 2\% to 5\%.

\subsection{Main Results}
\label{sec:main}

Here, we provide a more precise statement of our result \eqref{eq:heuristic}. We assume throughout that $\rf_s$ is a random forest built by training regression trees $T$ on subsamples of size $s$ out of $n$ drawn without replacement. Moreover, we always take the number of bootstrap replicates $B$ in \eqref{eq:rfm} to be large enough for Monte Carlo effects not to matter. \citet{mentch2014ensemble} study the impact of Monte Carlo noise on subsampled estimators; the Monte Carlo distribution of $\hVIJ$ is discussed by \citet{wager2014confidence}.

The following theorem summarizes our main contributions, which are described in more detail in Section \ref{sec:theor}. The consistency result \eqref{eq:intro_consitent} is comparable to theorems of \citet{meinshausen2006quantile}, \citet{biau2008consistency}, and \citet{biau2012analysis}. We are not aware, however, of any prior results resembling \eqref{eq:intro_gauss} or \eqref{eq:intro_ij}. We postpone our discussion of the regularity conditions on the base learner $T$ until Section \ref{sec:trees}.

\begin{theo}
\label{theo:intro}
Suppose that we have $n$ independent and identically distributed training examples $Z_i = \p{X_i, \, Y_i} \in [0, \, 1]^d \times [-M, \, M]$ for some $M > 0$. Suppose moreover that the distribution of the $X_i$ admits a density that is bounded away from both zero and infinity, that $\EE{Y^k \cond X = x}$ is Lipschitz-continuous in $x$ for $k = 1$ and $2$, and that we are interested in a test point $x$ at which $\Var{Y \cond X = x} > 0$.

Given this data-generating process, let $T$ be an honest, regular tree in the sense of Definitions \ref{defi:regular} and \ref{defi:honest} from Section \ref{sec:trees}, and let $\rf_{s\p{n}}$ be a random forest with base learner $T$ and a subsample size $s\p{n}$. Then, provided that the subsample size $s\p{n}$ satisfies
\begin{equation}
\label{eq:intro_cond}
{s\p{n}} \rightarrow \infty \; \eqand \;  \frac{s\p{n} \log\p{n}^{d}}{n} \rightarrow 0,
\end{equation}
the $B \rightarrow \infty$ random forest will be consistent:
\begin{equation}
\label{eq:intro_consitent}
\hy_n := \rf_{s\p{n}}\p{x; \, Z_1, \, ..., \, Z_n} \rightarrow_p \mu\p{x}
\end{equation}
as $n \rightarrow \infty$, where $\mu\p{x} = \EE{Y \cond X = x}$. Moreover, there exists a sequence $\sigma_n \rightarrow 0$ such that
\begin{equation}
\label{eq:intro_gauss}
\p{\hy_n -  \EE{\hy_n}} \big/ {\sigma_n}  \Rightarrow \nn\p{0, \, 1},
\end{equation}
and this sequence can consistently be estimated using the infinitesimal jackknife \eqref{eq:hvij}:
\begin{equation}
\label{eq:intro_ij}
\hVIJ\p{x; \, Z_1, \, ..., \, Z_n}  \big/ \sigma_n^2 \rightarrow_p 1.
\end{equation}
\end{theo}

Here, the choice of subsample size $s$ underlies a bias-variance trade-off. The variance of a random forest is governed by two factors: if the individual trees have variance $v$ and correlation $\rho$, the random forest itself will have variance $\rho v$ \citep{breiman2001random,hastie2009elements}. By making $s$ small we can reduce the overlap between different subsamples, thus decreasing $\rho$ and bringing down the variance of the ensemble. Conversely, when $s$ is large, the trees can grow deep and get very close to being unbiased.

\textbf{Remark: Bias.} It is often difficult to characterize the asymptotic bias of statistical predictors. Our case is no exception to this: Theorem \ref{theo:intro} guarantees that the bias $\EE{\hy} - \mu\p{x}$ goes to zero, but does not describe the scaled bias $\beta_n = \p{\EE{\hy} - \mu\p{x}} / \sigma_n(x)$. In general, we expect $\beta_n$ to converge to zero provided that the ratio $s\p{n} / n$ does not decay too fast. Establishing connections between the decay rate of $\beta_n$ and the smoothness of $\mu\p{x}$ would be an interesting avenue for further research.

\textbf{Addendum.} This project was in fact carried out in the revised paper by Wager and Athey, who derive a centered central limit theorem for random forests.

\subsection{Growing Trees}
\label{sec:trees}

In order for the results described in Theorem \ref{theo:intro} to hold, we of course need to put some constraints on the class of trees $T$ that can be used as base learners. We start by stating some regularity conditions below. The stronger and more interesting condition on $T$ is an ``honesty'' condition described in Definition \ref{defi:honest}. The trees below are grown on a subsample $\mathcal{S} \subseteq \{1, \, ..., \, n\}$ of the training examples.

\begin{defi}
\label{defi:regular}
A tree predictor grown by recursive partitioning is called regular if
\begin{enumerate}[\ \ \ (A)]
\item At each step of the training algorithm, the probability that the tree splits on variable $j$ is bounded below by some $\pi > 0$ for all $j = 1, \, ..., \, d$.
\item Each split leaves at least a fraction $\gamma > 0$ of the available training examples on each side of the split.
\item The trees are fully grown in the following sense. There is a set of prediction points $\mathcal{P} \subseteq \mathcal{S}$ of size $|\mathcal{P}| \geq |\mathcal{S}|/2$ such the tree makes predictions of the form
$ T\p{x} = Y_{i^*\p{x}}, $
where the index $i^*\p{x}$ corresponds to the sole $i \in \mathcal{P}$ such that $X_i$ and the test point $x$ are in the same leaf.
\end{enumerate}
\end{defi}

Here, conditions (A) and (B) are technical devices introduced by \citet{meinshausen2006quantile} to make sure that the predictions made by random forests become local as the trees get deep. Meanwhile, (C) is a theoretical convenience that lets us simplify the exposition. In practice, trees are sometimes grown to have terminal node size $k$ rather rather than 1 for regularization. In our setup, however, we already get regularization by drawing subsamples of size $s$ where $s \log(n)^d /n \rightarrow 0$ and the regularization effect from using larger leaf sizes is not as important.

Our second requirement is that $T$ be ``honest''. It is well known that trees act as adaptive nearest neighbors predictors; this idea is discussed in detail by \citet{lin2006random}. As stated more formally below, we say that a tree is honest if it does not re-use training labels $Y_i$ for both determining split-points of the tree and for making nearest-neighbor predictions. As far as we know, prior papers on the consistency of random forests \citep[e.g.,][]{biau2012analysis,biau2008consistency,meinshausen2006quantile} all require a similar condition.

\begin{defi}
\label{defi:honest}
A fully grown tree is called honest if, conditionally on $X_i$, the distribution of $Y_i$ does not depend on knowing whether $i$ is the selected index $i^*\p{x}$:
\begin{equation}
\label{eq:honest}
\law\p{Y_{i} \cond X_{i} = x, \, i^*(x) = i} \eqd  \law\p{Y_i \cond X_i = x}
\end{equation}
for all values of $x$.
\end{defi}

The simplest way to enforce the above condition is to divide the available training points $\mathcal{S}$ into a set of structure points $\mathcal{T}$ that are only used to pick tree splits and a set of prediction points $\mathcal{P}$ that are only used to make adaptive nearest neighbors predictions once we have already fixed the splits.\footnote{In order to satisfy condition (C), it may be necessary to introduce additional splits over the prediction set to make sure that each leaf has only one element. This can be done randomly.}

If we do not require \eqref{eq:honest}, it is easy to construct arbitrarily biased tree classifiers $T$. Suppose, for example, that $Y_i \in \{0, \, 1\}$ and that we are interested in a specific test point $x$. Then, a dishonest tree $T$ can always predict $\hy = 1$ if there exists even a single test example $(X_i, \, Y_i)$ such that $Y_i = 1$ and $X_i$ is a potential nearest neighbor (PNN) of $x$, i.e., there are no other test examples in the smallest axis-aligned rectangle containing both $x$ and $X_i$. As shown by \citet{lin2006random}, the expected number of PNNs of $x$ grows to infinity with $n$, and so in general, for large $n$, we can with high probability find a tree such that $T(x) = 1$.

\spacingset{\SPACESMALL}

\begin{figure}[t]
\centering
\includegraphics[width=\figw]{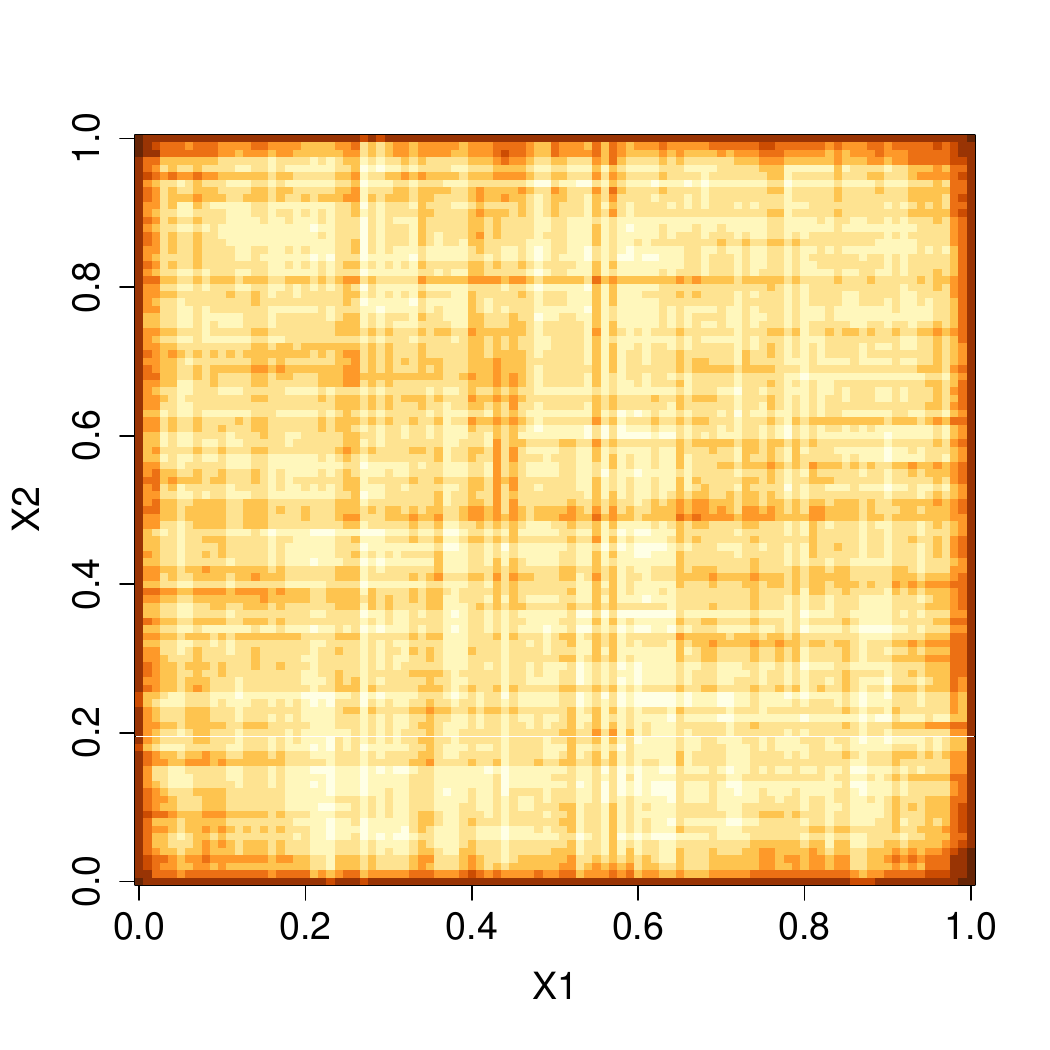}
\caption{Heat map of predictions made by a CART random forest with $n = 1,000,000$ training examples and a subsample size $s = 1,000$. The optimal prediction function should be constant, but the random forest systematically makes larger predictions near the edges and small ones in the middle of the domain.}
\label{fig:cart_bias}
\end{figure}

\spacingset{\SPACEBIG}

CART trees are of course not honest in the sense of Definition \ref{defi:honest}, because they use training labels both to choose splits and make predictions. And, as we show in Figure \ref{fig:cart_bias}, random forests trained with CART trees are not consistent either even for very simple problems. In this example, we drew $X$ uniformly from $[0, \, 1]^2$, and $Y \sim \operatorname{Bernoulli}\p{0.01}$ independently of $X$; thus, $\EE{Y \cond X = x} = 0.01$ everywhere. CART-based random forests, however, do not consistently learn this constant function, but rather seem to warp up near the edges and down in the middle. At $x = (0, \, 0)$, the tree was predicting $\hy = 3.5\%$ instead of $1\%$. It appears that the CART trees are trying to aggressively separate the points with $Y_i = 1$ from the rest of the data, and in doing so push the neighborhoods surrounding those points towards the edge of the domain of $X$.

That being said, the bias of CART trees seems to be subtle enough that it does not affect the performance of random forests in most situations. Similarly, in our simulation experiments presented in Section \ref{sec:sim}, we find that the infinitesimal jackknife $\hVIJ$ works well for estimating the variance of random forests with CART base learners. Thus, it seems reasonable to assume that our main results from Theorem \ref{theo:intro} still provide useful insight in understanding the behavior of CART random forests as implemented in, e.g., the popular \texttt{R} library \texttt{randomForest} \citep{liaw2002classification}.

Seeing whether it is possible to improve the practical performance of CART trees by making them honest and unbiased seems like a promising avenue for further research. In a recent advance, \citet{denil2014narrowing} proposed a form of random forests that satisfy Definition \ref{defi:honest} while matching the empirical performance of CART-based random forests on several datasets.

\section{Theoretical Development}
\label{sec:theor}

The ideas used in our proof go back to techniques developed by \citet{hoeffding1948class} and \citet{hajek1968asymptotic} to establish the asymptotic normality of classical statistical estimators such as $U$-statistics. We begin by briefly reviewing their results to give some context to our proof. Given a predictor $T$ and independent training examples $Z_1$, ..., $Z_n$, the H\'ajek projection of $T$ is defined as
\begin{equation}
\label{eq:hajek}
\proj{T} = \EE{T} + \sum_{i = 1}^n \p{\EE{T \cond Z_i} - \EE{T}}.
\end{equation}
In other words, the H\'ajek projection of $T$ captures the first-order effects in $T$.
This projection has the properties we would expect: in particular $\Var{\proj{T}} \leq \Var{T}$, and if
\begin{equation}
\label{eq:hajek_l2}
\limn {\Var{\proj{T}}} \big/ {\Var{T}} = 1, \text{ then } \limn {\EE{\Norm{\proj{T} - T}_2^2}} \big/ {\Var{T}} = 0. 
\end{equation}
Since the H\'ajek projection $\proj{T}$ is a sum of independent random variables, we should expect it to be asymptotically normal under all but pathological conditions. Thus whenever the ratio of the variance of $\proj{T}$ to that of $T$ tends to 1, the theory of H\'ajek projections almost automatically guarantees that $T$ will be asymptotically normal.\footnote{The moments defined in \eqref{eq:hajek} depend on the data-generating process for the $Z_i$, and so cannot be observed in practice. Thus, the H\'ajek projection is mostly useful as an abstract theoretical tool. For a review of classical projection arguments, see Chapter 11 of \citet{van2000asymptotic}.}

If $T$ is decision tree, however, the condition from \eqref{eq:hajek_l2} does not apply, and we cannot use the classical theory of H\'ajek projections directly. Our analysis is centered around a weaker form of this condition, which we call $\alpha$-incrementality. With our definition, predictors $T$ to which we can apply the argument \eqref{eq:hajek_l2} directly are 1-incremental.

\begin{defi}
\label{defi:incr}
The predictor $T$ is $\alpha(s)$-incremental at $x$ if
$$ {\Var{\proj{T}\p{Z_1,\, ..., \, Z_s}}} \big / {\Var{T\p{Z_1,\, ..., \, Z_s}}} \gtrsim \alpha(s), $$
where $\proj{T}$ is the H\'ajek projection of $T$ \eqref{eq:hajek}. In our notation,
$$f(s) \gtrsim g(s) \text{ means that } \liminf_{s \rightarrow \infty} {f(s)} \big/ {g(s)} \geq 1. $$
\end{defi}

Our argument proceeds in two steps. First we establish lower bounds for the incrementality of regression trees in Section \ref{sec:incremental}. Then, in Section \ref{sec:hajek} we show how we can turn weakly incremental predictors $T$ into 1-incremental ensembles by subsampling (Lemma \ref{lemm:hajek}), thus bringing us back into the realm of classical theory. We also establish the consistency of the infinitesimal jackknife for random forests.

Our analysis of regression trees is motivated by the ``potential nearest neighbors'' model for random forests introduced by \citet{lin2006random}. Meanwhile, the key technical device used in Section \ref{sec:hajek} is the ANOVA decomposition of \citet{efron1981jackknife}. The discussion of the infinitesimal jackknife for random forest builds on results of \citet{efron2013estimation} and \citet{wager2014confidence}.

\subsection{Regression Trees and Incremental Predictors}
\label{sec:incremental}

Analyzing specific greedy tree models such as CART trees can be challenging. We thus follow the lead of \citet{lin2006random}, and analyze a more general class of predictors---potential nearest neighbors predictors---that operate by doing a nearest-neighbor search over rectangles.

\begin{defi}
Consider a set of points $X_1, \, ..., \, X_s \in \RR^d$ and a fixed $x \in \RR^d$. A point $X_i$ is a potential nearest neighbor (PNN) of $x$ if the smallest axis-aligned hyperrectangle with vertices $x$ and $X_i$ contains no other points $X_j$. A predictor $T$ is a PNN predictor if, given a training set $\p{X_1, \, Y_1}, \, ..., \, \p{X_s, \, Y_s} \in \RR^d \times \yy$ and a test point $x \in \RR^d$, $T$ always outputs $Y_i$ corresponding to a PNN $X_i$ of $x$.
\end{defi}

In other words, given a test point $x$, a PNN classifier is allowed to predict $Y_i$ only if there exists a rectangle containing $x$, $X_i$, and no other test points. All standard decision trees result in PNN predictors.\footnote{It is straight-forward to extend our results to the case where trees have leaves with at most $k$ elements in each leaf, e.g., $k = 5$ for \texttt{randomForest} \citep{liaw2002classification} in regression mode; see \citet{lin2006random}.}

\begin{prop}[\citet{lin2006random}]
Any decision tree $T$ that makes axis-aligned splits and has leaves of size 1 is a PNN predictor. In particular, the base learners originally used by \citet{breiman2001random}, namely fully grown CART trees \citep{breiman1984classification}, are PNN predictors.
\end{prop}

Predictions made by PNNs can always be written as
\begin{equation}
\label{eq:pnn_sum}
T\p{x; \, Z} = \sum_{i = 1}^s S_i Y_i,
\end{equation}
where $S_i$ is a selection variable that takes the value 1 for the selected index $i$ and 0 for all other indices. An important property of PNN predictors is that we can often get a good idea about whether $S_i$ can possibly be 1 even if we only get to see $Z_i$; more formally, as we show below, the quantity $n\Var{S_1\cond Z_1}$ cannot get too small. Establishing this fact is a key step in showing that PNNs are incremental.

\begin{lemm}
\label{lemm:pnn}
Suppose that the features $X$ are independently and identically distributed on $[0, \, 1]^d$ with a density $f$ that is bounded away from zero and infinity, and let $T$ be any PNN predictor. Then, there is a constant $C_f$ depending only of $f$ and $d$ such that
\begin{equation}
\label{eq:pnn_bound}
s\Var{\EE{S_1 \cond Z_1}} \gtrsim C_f \big/ \log\p{s}^{d},
\end{equation}
where $S_i$ is defined as in \eqref{eq:pnn_sum}.
When $f$ is uniform over $[0, \, 1]^d$, the bound holds with $C_f = 2^{d + 1} / \p{d - 1}!$.
\end{lemm}

Thanks to this result, we are now ready to show that all honest regular trees are incremental.  The proof of the following theorem makes use of an important technical lemma from \citet{meinshausen2006quantile}.
Finally, we note that the result from Theorem \ref{theo:pnn} can easily be extended to trees as described in part (C) of Definition \ref{defi:regular} and recommended by \citet{denil2014narrowing}, where the tree is only a potential nearest neighbors classifier over a set of prediction points $\mathcal{P} \subseteq \mathcal{S}$ of size $\mathcal{P} \geq s/2$.

\begin{theo}
\label{theo:pnn}
Suppose that the conditions of Lemma \ref{lemm:pnn} hold and that $T$ is an honest regular tree in the sense of Definitions \ref{defi:regular} and \ref{defi:honest}. Suppose moreover that the moments $\EE{Y\cond X = x}$ and $\EE{Y^2 \cond X = x}$ are both finite and Lipschitz continuous functions of $x$, and that $\EE{Y^2 \cond X = x}$ is uniformly bounded for all $x \in [0, \, 1]^d$. Finally, suppose that $\Var{Y \cond X = x} > 0$. Then $T$ is $\alpha\p{s}$-incremental with
$$ \alpha\p{s} = {C_f} \big / {\log\p{s}^{d}}, $$
where $C_f$ is the constant from Lemma \ref{lemm:pnn}.
\end{theo}

\subsection{Random Forests with Incremental Base Learners}
\label{sec:hajek}

In the previous section, we showed that decision trees are $\alpha$-incremental, in that the H\'ajek projection $\proj{T}$ of $T$ contains at least some information about $T$. In this section, we show that randomly subsampling $\alpha$-incremental predictors makes them 1-incremental; this then lets us proceed with a classical statistical analysis. The following lemma, which flows directly from the ANOVA decomposition of \citet{efron1981jackknife}, provides a first motivating result for our analysis.

\begin{lemm}
\label{lemm:hajek}
Let $\rf_s$ be a random forest with base learner $T$ as defined in \eqref{eq:rfm}, and let $\proj{\rf_s}$ be the H\'ajek projection of $\rf_s$ \eqref{eq:hajek}. Then
$$\EE{\p{\rf_s\p{x; \, Z_1, \, ..., \, Z_n}  - \wrf_s\p{x; \, Z_1, \, ..., \, Z_n} }^2} \leq \p{\frac{s\p{n}}{n}}^2 \, \Var{T} $$
whenever the variance $\Var{T}$ of the base learner is finite.
\end{lemm}

This technical result paired with Theorem \ref{theo:pnn} leads to our first main result about random forests.

\begin{theo}
\label{theo:gauss}
Let $\rf_{s\p{n}}$ be a random forest with base learner $T$ trained according the conditions of Theorem \ref{theo:pnn}, with $Y$ restricted to a bounded interval $Y \in [-M, \, M]$. Suppose, moreover, that the subsample size $s\p{n}$ satisfies
$$ \limn s\p{n} = \infty \eqand \limn {s\p{n} \log\p{n}^{d}} \big/ {n} = 0. $$
Then, the random forest is consistent.
Moreover, there exists a sequence $\sigma_n \rightarrow 0$ such that
\begin{equation}
\label{eq:gauss}
\frac{1}{\sigma_n} \p{\rf_{s\p{n}}\p{x; \, Z_1, \, ..., \, Z_n} - \EE{\rf_{s\p{n}}\p{x; \, Z_1,\, ...,\, Z_n}}} \Rightarrow \nn\p{0, \, 1},
\end{equation}
where $\nn\p{0, \, 1}$ is the standard normal density.
\end{theo}

Theorem \ref{theo:gauss} by itself is just an abstract characterization theorem about random forests. As we show below, however, it is possible to accurately estimate the noise-level $\sigma^2(\hy)$ of a random forest using the infinitesimal jackknife for random forests \citep{efron2013estimation,wager2014confidence}. Thus, our theory allows us to do statistical inference about random forest predictions.

\begin{theo}
\label{theo:ij}
Let $\hVIJ\p{x;, \, Z_1,\, ...,\, Z_n}$ be the infinitesimal jackknife for random forests \citep{efron2013estimation} as defined in \eqref{eq:hvij}. Then, under the conditions of Theorem \ref{theo:gauss},
$$ \hVIJ\p{x; \, Z_1,\, ...,\, Z_n} \big/ \sigma_n^2 \rightarrow_p 1. $$
\end{theo}

\section{Simulation Study}
\label{sec:sim}

In this section we test the result from Theorem \ref{theo:ij}, namely that the infinitesimal jackknife estimate of variance $\hVIJ$ is a good predictor of the true variance $\sigma^2(\hy)$ of random forest predictions. The simulations all used a subsample size of $s = \lfloor n^{0.7} \rfloor$ and $B = 5 \, n$ bootstrap replicates; otherwise, we used default settings for the \texttt{randomForest} package in \texttt{R} \citep{liaw2002classification}.

As emphasized by \citet{wager2014confidence}, when we only use $B = \oo(n)$ bootstrap replicates, $\hVIJ$ can suffer from considerable Monte Carlo bias. To alleviate this problem, they propose a Monte Carlo bias correction for $\hVIJ$. In the case of subsampling with subsample size $s$ and $B$ bootstrap replicates, this bias correction is
\begin{align}
\label{eq:bias_corr}
\hVIJ^B &= \sum_{i = 1}^n C_i^2 - \frac{s \p{n - s}}{n} \, \frac{\hv}{B}, \where \\
&C_i = \frac{1}{B} \sum_{b = 1}^B \p{N_{bi}^* - s/n} \, \p{T^*_b -\bar{T}^*} \eqand \\
&\hv = \frac{1}{B} \sum_{b = 1}^B \p{T^*_b - \bar{T}^*}^2.
\end{align}
Here, the first term of $\hVIJ^B$ is a plug-in formula where $C_i$ acts as a Monte Carlo estimate for $\Cov[*]{N_i^*, \, T^*_n}$, and the second term is a bias correction.

\spacingset{\SPACESMALL}

\begin{table}[t]
\begin{center}
\begin{tabular}{r|rr||ccc||ccc||}
& & &  \multicolumn{3}{c||}{Relative} & \multicolumn{3}{c||}{Absolute} \\
 Distr. & d & n & Bias$^2$ & Variance & MSE & Bias$^2$ & Variance & MSE \\ 
   \hline
\hline
Cosine & 2 & 200 & 0.13 & 0.08 & 0.21 & 1.23E-03 & 7.76E-04 & 2.00E-03 \\ 
  Cosine & 2 & 1000 & 0.04 & 0.08 & 0.11 & 3.97E-05 & 8.52E-05 & 1.24E-04 \\ 
  Cosine & 2 & 5000 & 0.03 & 0.05 & 0.08 & 6.44E-06 & 1.17E-05 & 1.80E-05 \\ 
   \hline
Cosine & 10 & 200 & 0.10 & 0.07 & 0.17 & 7.60E-04 & 5.63E-04 & 1.32E-03 \\ 
  Cosine & 10 & 1000 & 0.10 & 0.04 & 0.14 & 1.31E-04 & 5.06E-05 & 1.81E-04 \\ 
  Cosine & 10 & 5000 & 0.04 & 0.03 & 0.07 & 4.42E-06 & 3.28E-06 & 7.66E-06 \\
     \hline
\hline
AND & 20 & 200 & 0.08 & 0.09 & 0.17 & 6.92E-03 & 7.35E-03 & 1.42E-02 \\ 
  AND & 20 & 1000 & 0.22 & 0.07 & 0.29 & 2.06E-03 & 6.54E-04 & 2.71E-03 \\ 
  AND & 20 & 5000 & 0.65 & 0.27 & 0.92 & 3.01E-04 & 1.26E-04 & 4.26E-04 \\ 
   \hline
AND & 100 & 200 & 0.07 & 0.18 & 0.24 & 2.15E-03 & 5.69E-03 & 7.78E-03 \\ 
  AND & 100 & 1000 & 0.07 & 0.11 & 0.18 & 3.27E-04 & 5.00E-04 & 8.23E-04 \\ 
  AND & 100 & 5000 & 0.07 & 0.25 & 0.33 & 4.40E-05 & 1.49E-04 & 1.92E-04 \\  
   \hline
\hline
XOR & 5 & 200 & 0.10 & 0.07 & 0.17 & 7.08E-03 & 4.51E-03 & 1.16E-02 \\ 
  XOR & 5 & 1000 & 0.05 & 0.04 & 0.10 & 4.66E-04 & 3.95E-04 & 8.56E-04 \\ 
  XOR & 5 & 5000 & 0.06 & 0.02 & 0.08 & 5.38E-05 & 2.38E-05 & 7.73E-05 \\ 
   \hline
XOR & 20 & 200 & 0.12 & 0.07 & 0.19 & 7.58E-03 & 3.97E-03 & 1.15E-02 \\ 
  XOR & 20 & 1000 & 0.06 & 0.04 & 0.09 & 5.38E-04 & 3.27E-04 & 8.62E-04 \\ 
  XOR & 20 & 5000 & 0.04 & 0.02 & 0.06 & 4.53E-05 & 2.56E-05 & 7.07E-05 \\ 
   \hline
\hline
\end{tabular}
\vspace{2mm}
\caption{Performance of the infinitesimal jackknife for random forests on synthetic distributions. The ``absolute'' metrics describe the accuracy of $\hVIJ$, while the ``relative'' metrics describe the accuracy of $\hVIJ/\bar{\sigma}^2$, where $\bar{\sigma}^2$ is the average of $\sigma^2(\hy)$ over the test set. All examples have a subsample size of $s = \lfloor n^{0.7} \rfloor$ and use $B = 5n$ bootstrap replicates.}
\label{tab:rf_simu}
\end{center}
\end{table}

\spacingset{\SPACEBIG}

Table \ref{tab:rf_simu} shows the performance of $\hVIJ$ on several synthetic distributions described at the end of this section. To produce this table, we first drew $K = 100$ random test points $\{x^{(k)}\}_{k = 1}^{K}$ from the data-generating distribution. We then constructed $R = 100$ random training sets $\{\bZ^{(r)}\}_{r = 1}^{R}$, and evaluated both the prediction $\rf_s\p{x^{(k)}; \, \bZ^{(r)}}$ and the variance estimate $\hVIJ\p{x^{(k)}; \, \bZ^{(r)}}$ for each test point $x^{(k)}$. The numbers shown in Table \ref{tab:rf_simu} are averaged over the test points $k$:
\begin{align*}
&\text{Bias}^2 = \frac{1}{k} \sum_{k = 1}^K \p{\frac{1}{R} \sum_{r = 1}^R \hVIJ^B\p{x^{(k)}; \, \bZ^{(r)}} - \Var[r]{\rf_s\p{x^{(k)}; \, \bZ^{(r)}}} }^2, \\
&\text{Var} =\frac{1}{k} \sum_{k = 1}^K \frac{1}{R - 1} \sum_{r = 1}^R \p{\hVIJ^B\p{x^{(k)}; \, \bZ^{(r)}} - \frac{1}{R} \sum_{r = 1}^R \hVIJ^B\p{x^{(k)}; \, \bZ^{(r)}}}^2,
\end{align*}
and $\text{MSE} = \text{Bias}^2 + \text{Var}$.
In the ``absolute'' columns we show the raw metrics as is, while for the ``relative'' columns we divided the metrics by the average of $\sigma^2(\hy)^2$ over the test points.

Rather encouragingly, we see that $\hVIJ$ is overall quite accurate, and gets more accurate as $n$ gets larger. For the ``Cosine'' and ``XOR'' distributions, the relative MSE also decays with $n$, as predicted by Theorem \ref{theo:ij}. The ``AND'' distribution appears to have been the most difficult distribution: although the error of $\hVIJ$ decays with $n$ as seen in the ``Absolute'' columns, we have not yet entered the regime where it decays faster than $\sigma^2(\hy_n)$. The fact that the ``AND'' distribution would be the most difficult one is not surprising, as it has the highest dimensionality $d$, and the conditional mean function is far from being Lipschitz continuous.

We also tested $\hVIJ$ on some parametric bootstrap simulations based on classic data sets from the UCI machine learning repository \citep{UCI}; results are shown in Table \ref{tab:rf_simu_uci}. The ``forest fires'' and ``housing'' data sets are originally due to \citet{cortez2007data} and \citet{harrison1978hedonic}; with the forest fires data set, we predicted log-area of the fire.
For these experiments, we divided the original data set into a test set and a training set, and then used the training set to construct a parametric bootstrap distribution. From there, we proceeded as with the synthetic distributions from Table \ref{tab:rf_simu}. The infinitesimal performed quite well here, despite the small sample size.

\spacingset{\SPACESMALL}

\begin{table}[t]
\begin{center}
\begin{tabular}{r|rr||ccc||ccc||}
& & &  \multicolumn{3}{c||}{Relative} & \multicolumn{3}{c||}{Absolute} \\
 Distr. & d & n & Bias$^2$ & Variance & MSE & Bias$^2$ & Variance & MSE \\ 
  \hline 
 	Auto & 7 & 314 & 0.09 & 0.10 & 0.19 & 5.30E-03 & 5.90E-03 & 1.11E-02 \\ 
   Fires & 12 & 344 & 0.10 & 0.15 & 0.25 & 4.73E-04 & 6.90E-04 & 1.16E-03 \\ 
  Housing & 13 & 337 & 0.09 & 0.14 & 0.23 & 2.46E-02 & 4.05E-02 & 6.46E-02 \\ 
   \hline
\end{tabular}
\vspace{2mm}
\caption{Performance of the infinitesimal jackknife for random forests on the auto, forest fires and Boston housing data sets from the UCI repository \citep{UCI}.  The ``absolute'' metrics describe the accuracy of $\hVIJ$, while the ``relative'' metrics describe the accuracy of $\hVIJ/\bar{\sigma}^2$, where $\bar{\sigma}^2$ is the average of $\sigma^2(\hy)$ over the test set. All examples have a subsample size of $s = \lfloor n^{0.7} \rfloor$ and use $B = 5n$ bootstrap replicates.}
\label{tab:rf_simu_uci}
\end{center}
\end{table}

\spacingset{\SPACEBIG}

\subsection{Synthetic Data-Generating Distributions}

To generate the data used in Table \ref{tab:rf_simu}, we first drew features $X \sim U\p{[0, \, 1]^p}$ and then generated the labels $Y$ using the following rules:
\begin{align*}
&\text{Cosine: } Y = 3 \cdot \cos\left(\pi \cdot (X_1 + X_2)\right) + \varepsilon \\
&\text{XOR: } Y = 5 \cdot \left[\texttt{XOR}\left(X_1 > 0.6, \, X_2 > 0.6\right) + \texttt{XOR}\left(X_3 > 0.6, \, X_4 > 0.6\right)\right] + \varepsilon  \\
&\text{AND: } Y = 10 \cdot \texttt{AND}\left(X_1 > 0.3, \, X_2 > 0.3, \, X_3 > 0.3, \, X_4 > 0.3 \right) + \varepsilon.
\end{align*}
Here, $\varepsilon$ is standard Gaussian noise, while $\texttt{XOR}$ and $\texttt{AND}$ are treated as 0/1-valued functions. Similar distributions were also used for simulation experiments in \citep{wager2014confidence}.

\section{Conclusion}

Following Breiman's lead, random forests have mostly been studied in terms of black-box metrics like cross-validation error or test-set error. Although measuring test-set error is of course important, a pure black-box evaluation can lack much of the richness of a more careful statistical analysis. It cannot answer questions like: ``Are there some points at which the random forest is more stable than others?,'' or ``How much might the predictions change if we gathered more training data?''

In this paper, we studied a random forest model based on subsampling with a subsample size satisfying $s/n = o(\log(n)^{-d})$. We established results about the asymptotic normality of random forest predictions, and showed how their asymptotic variance can be accurately estimated from data. Thus, our results help open the door to a more complete and nuanced statistical understanding of random forest models.

\spacingset{\SPACESMALL}

\bibliographystyle{plainnat}
\bibliography{references}

\end{document}